\documentclass[11pt]{article}
\usepackage{amsmath} \usepackage{amsthm} \usepackage{amsfonts}
\usepackage{graphicx}
\usepackage{color}

\newtheorem{theorem}{Theorem} 
\newtheorem{lemma}{Lemma} 
\newtheorem{cor}{Corollary}

\newcommand{\indic}[1]{\mathbf{1}_{\{#1\}}}

\def\Nn{\mathbb{N}}

\newcommand{\tendsto}{\to}
\newcommand{\Cc}{\mathcal{C}}
\newcommand{\ds}{\displaystyle}
\newcommand{\cvlaw}{\Rightarrow}
\newcommand{\var}{\text{var}\,}
\newcommand{\cov}{\text{cov}\,}

\newcommand{\keywords}{\textbf{Keywords}\ }

\def\erm{\footnotesize}
\def\mn{\medskip\noindent}
\def\bn{\bigskip\noindent}
\def\beq{\begin{equation}}
\def\eeq{\end{equation}}
\def\beqa{\begin{eqnarray}}
\def\eeqa{\end{eqnarray}}
\def\beqax{\begin{eqnarray*}}
\def\eeqax{\end{eqnarray*}}
\def\sqz{\kern -0.2em}

\begin{document}

\title{A phase transition in the random transposition random walk}

\author{Nathana\"el Berestycki$^{1,2}$ and Rick Durrett$^2$}
\maketitle

\begin{abstract}
Our work is motivated by Bourque and Pevzner's (2002) simulation
study of the effectiveness of the parsimony method in studying
genome rearrangement, and leads to a surprising result about the
random transposition walk on the group of permutations on $n$
elements. Consider this walk in continuous time starting at the
identity and let $D_t$ be the minimum number of transpositions
needed to go back to the identity from the location at time $t$.
$D_t$ undergoes a phase transition: the distance $D_{cn/2} \sim
u(c)n$, where $u$ is an explicit function satisfying $u(c)=c/2$
for $c \le 1$ and $u(c)<c/2$ for $c>1$. In addition, we describe
the fluctuations of $D_{cn/2}$ about its mean in each of the three
regimes (subcritical, critical and supercritical). The techniques
used involve viewing the cycles in the random permutation as a
coagulation-fragmentation process and relating the behavior to the
Erd\H{o}s-Renyi random graph model.
\end{abstract}

\vfill\noindent \keywords{random transposition, random graphs,
phase transition, coagulation-fragmentation, genome rearrangement,
parsimony method}

\bn 1. Ecole Normale Sup\'erieure, D\'epartement de
Math\'ematiques et Applications, 45, rue d'Ulm F-75005 Paris,
France

\bn 2. Department of Mathematics, Malott Hall, Cornell University,
Ithaca, NY 14853, U.S.A. Both authors are partially supported by a
joint NSF-NIGMS grant DMS-0201037, and would like to thank David
Aldous for taking the time to answer a number of questions as this
paper was being written. \clearpage

\section{General motivation}
The relationship between the orders of genes in two species can be
described by a signed permutation. For example the relationship
between the human and mouse $X$ chromosomes may be encoded as (see
Pevzner and Tesler (2003))
$$
1 \quad {-7} \quad 6 \quad {-10} \quad 9 \quad {-8} \quad 2 \quad
{-11} \quad {-3} \quad 5 \quad 4
$$
In words the two $X$ chromosomes can be partitioned into 11
segments. The first segment of the mouse $X$ chromosome is the
same as that of humans, the second segment of mouse is the 7th
human segment with its orientation reversed, etc. The parsimony
approach to estimation of evolutionary changes of the $X$
chromosome between human and mouse is to ask: what is the minimum
number of reversals (i.e., moves that reverse the order of a
segment and therefore change its {\it sign}) needed to transform
the arrangement above back into $1, \ldots, 11$? In other words,
what is the (reversal) distance between the human and mouse X
chromosomes ?

Hannehalli and Pevzner (1995) developed a polynomial algorithm for
answering this question. The first step in preparing to use the
Hannehalli-Pevzner algorithm is to double the markers. When
segment $i$ is doubled we replace it by two consecutive numbers
$2i-1$ and $2i$, e.g., 6 becomes 11 and 12. A reversed segment
$-i$ is replaced by $2i$ and $2i-1$, for example, $-7$ is replaced
by 14 and 13. The doubled markers use up the integers 1 to 22. To
these numbers we add a 0 at the front and a 23 at the end. Using
commas to separate the ends of the markers we can write the two
genomes as follows:
\begin{align}
\hbox{mouse} \quad & 0, 1\, 2, 14 \, 13, 11 \, 12, 20 \, 19, 17 \,
18, 16 \, 15, 3 \, 4, 22 \, 21,
6 \, 5, 9 \, 10, 7 \, 8, 23 \nonumber \\
\hbox{human} \quad & 0, 1 \, 2, 3 \, 4, 5 \, 6, 7 \, 8, 9\, 10, 11
\, 12, 13 \, 14, 15 \, 16, 17 \, 18, 19 \, 20, 21 \, 22,
23\nonumber
\end{align}
The next step is to construct the breakpoint graph (see Figure 1)
that results when the commas are replaced by edges that connect
vertices with the corresponding numbers. In the picture we have
written the vertices in their order in the mouse genome. Commas in
the mouse order become thick lines (black edges), while those in
the human genome are thin lines (gray edges).

Each vertex has one black and one gray edge, so the connected
components of the graph are easy to find: start with a vertex and
follow the connections in either direction until you come back to
where you start. In this example there are five components:
\begin{align}
& \kern0.5em 0-1-0 \qquad 2-14-15-3-2 \qquad 4-22-23-8-9-5-4  \cr
&19-17-16-18-19 \qquad 13-11-10-7-6-21-20-12-13\nonumber
\end{align}

To compute a lower bound for the distance, we take the number of
commas seen when we write out one genome. In this example that is
12. In general, it is 1 plus the number of markers. We then
subtract the number of components in the breakpoint graph. In this
example that is 5, so the result is 7. This is a lower bound on
the distance, since any reversal can at most reduce this quantity
by 1, and it is 0 when the two genomes are the same. We can verify
that 7 is the minimum distance by constructing a sequence of 7
moves that transforms the mouse $X$ chromosome into the human
order. There are thousands of solutions, so we leave this as an
exercise for the reader. Here are some hints: (i) To do this it
suffices, at each step, to choose a reversal that increases the
number of cycles by 1. (ii) This never occurs if the two chosen
black edges are in different cycles. (iii) If the two black edges
are in the same cycle and are $(a,b)$ and $(c,d)$ as we read from
left to right, this will occur unless in the cycle minus these two
edges $a$ is connected to $d$ and $b$ to $c$, in which case the
number of cycles will not change. For example, in the graph in
Figure 1 a reversal that breaks black edges 19-17 and 18-16 will
increase the number of cycles but the one that breaks 2-14 and
15-3 will not.

In general, the distance between genomes can be larger than the
lower bound from the breakpoint graph. There can be obstructions
called {\it hurdles} that can prevent us from decreasing the
distance, and hurdles can be intertwined in a {\it fortress of
hurdles} that takes an extra move to break. See Hannehalli and
Pevzner (1995). In symbols, if $\pi$ is the signed permutation that
represents the relative order and orientation of segments in the
two genomes, then
$$
d(\pi) = n+1 - c(\pi) + h(\pi) + f(\pi)
$$
where $d(\pi)$ is the distance from the identity, $n$ is the
number of markers, $c(\pi)$ is the number of components in the
breakpoint graph, $h(\pi)$ is the number of hurdles, and $f(\pi)$
is the indicator of the event $\pi$ is a fortress of hurdles. See
Section 5.2 of Durrett (2002) or Chapter 10 of Pevzner (2000) for
more details.

Although $d_0(\pi) = n+1 - c(\pi)$ is only a lower bound on the
distance, it is the right answer in most biological examples.
Bafna and Pevzner (1995) consider 11 comparisons of mitochondrial
and chloroplast genomes and found that this lower bound gave the
right answer in all cases. This pattern has continued in more
recent work, see York, Durrett, and Nielsen (2002), and Durrett,
Nielsen, and York (2003). The simulations in Figure 2 will give
more evidence that $d_0(\pi)$ and $d(\pi)$ are close in many
cases.

To motivate our main question, we will introduce a second data
set. Ranz, Casals, and Ruiz (2001) located 79 genes on chromosome
2 of {\it D. repleta} and on chromosome arm 3R of {\it D.
melanogaster}. If we number the genes according to their order in
{\it D. repleta} then their order in {\it D. melanogaster} is
given in Table 1. This time we do not know the orientation of the
segments, but that is not a serious problem. Using simulated
annealing, one can easily find an assignment of signs that
minimizes the distance, which in this case is 54. Given the large
number of rearrangements relative to the number of markers, we
should ask: when is the parsimony estimate reliable?

Bourque and Pevzner (2002) approached this question by taking 100
markers in order, performing $k$ randomly chosen reversals to get
a permutation $\pi_k$, computing the minimum number of reversals
needed to return to the identity, $d(\pi_k)$, and then plotting
the average value of $d(\pi_k) - k \le 0$ for 100 simulations.
They concluded, based on their simulations, that the parsimony
distance for $n$ markers was a good estimate as long as the number
of reversals performed was at most $0.4n$. In Figure 2 we have
given $-1$ times their data. We have also repeated their
experiment for the approximate distance $d_0(\pi) = n+1 - c(\pi)$
and plotted the average value of $k-d_0(\pi_k) \ge 0$ for 10,000
replications. Our curve is less random, but close to data of
Bourque and Pevzner (2000). The smooth curve gives result of
Theorem 3 for the limiting behavior of $(tn - d_0(\pi_{tn}))/n$
(as a function of $t$).

The biological question concerns the random reversal walk.
However, it is also interesting to consider the analogous problem
for random transpositions. In that case the distance from the
identity can be easily computed: it is the number of markers $n$
minus the number of cycles in the permutation. For an example,
consider the following permutation of 14 objects written in its
cyclic decomposition:
$$
(1\, 7\, 4) \: (2) \; (3 \, 12) \: (5 \, 13 \, 9 \, 11 \, 6) \; (8
\, 10 \, 14)
$$
which indicates that $1 \to 7$, $7 \to 4$, $4\to 1$, $2 \to 2$,
$3\to 12$, $12 \to 3$, etc. There are 5 cycles so the distance
from the identity is 9. If we perform a transposition that
includes markers from two different cycles (e.g., 7 and 9) the two
cycles merge into 1, while if we pick two in the same cycle (e.g.,
13 and 11) it splits into two.

The situation is similar but slightly more complicated for
reversals. There a reversal that involves edges in two different
components merges them into 1, but a reversal that involves two
edges of the same cycle may or may not increase the number of
cycles. One can attempt to couple the components of the breakpoint
graph for random reversals on $n-1$ markers and the cycles of
random transposition of $n$ markers as follows: number the edges
between markers in the reversal chain (including the ends 0 and
$n$); when markers $i$ and $j$ are transposed, do the inversion of
edges numbered $i$ and $j$. The result of the coupled simulation
is given in Figure 2. As expected time minus distance is smaller
for reversals but the qualitative behavior is similar. Thus, we
will begin by considering the biologically less relevant case of
random transpositions, and ask a question that in terms of the
rate 1 continuous time random walk on the symmetric group is: how
far from the identity are we at time $cn/2$? We will see later
that parts of the answer can be extended to the reversal random
walk.

\section{The coagulation-fragmentation process and the random graph process}

Let $(\sigma_t, t \ge 0)$ be the continuous-time random walk on
the group of permutations, starting at the identity, in which, at
times of a rate one Poisson process, we perform a transposition of
two elements chosen uniformly at random, with replacement, from
$\{1,\dots, n\}$. Choosing with replacement causes the chain to do
nothing with probability $1/n$, but makes some of the calculations
a little nicer. If we think of the permutation $\sigma$ as being
represented by numbered balls sitting on numbered locations with
ball $\sigma(k)$ sitting at $k$, then transposition of $i$ and
$j$, $\rho_{i,j}$, can be implemented in two ways. We can exchange
the balls at $i$ and $j$ or the balls numbered $i$ and $j$.
Algebraically these correspond to $\rho_{i,j}\sigma$ and
$\sigma\rho_{i,j}$. Since $(\sigma \rho_{i,j})^{-1} =
\rho_{i,j}\sigma^{-1}$ and the partition of $\{1, \dots, n\}$
induced by the cycle decompositions of $\sigma$ and $\sigma^{-1}$
are equal, the results are the same for either random walk.

\medskip Define the distance to the identity $D_t$ to be the
minimum number of transpositions one needs to perform on
$\sigma_t$ to go back to the identity element. A different way of
looking at $D_t$ is the following. $(\sigma_t, t\ge 0)$ can be
viewed as a random walk on a graph $G$, where $G$ is the Cayley
graph of the symmetric group for the set of generators given by
the set of all transpositions. Using this language, we see that
$D_t$ is nothing but the graph distance from $\sigma_t$ to the
origin, the identity element.

\medskip It is clear that if $N_t$ is the number of transpositions
distinct from the identity performed up to time $t$ (a Poisson
random variable with mean $t[1-2/(n-1)]$), then $D_t \le N_t$. As
mentioned earlier $D_t$ is given by $D_t = n-|\sigma_t|$, where
$|\sigma_t|$ is the number of cycles in the cycle decomposition of
$\sigma_t$. This formula allows us to turn any question about
$D_t$ into a question about $|\sigma_t|$. The key to studying
$|\sigma_t|$ is that the cycles evolve according to the dynamics
of a coagulation-fragmentation process. When a transposition
$\rho_{i,j}$ occurs, if $i$ and $j$ belong to two different cycles
then the cycles merge. On the contrary, if they belong to the same
cycle, this cycle is split into two cycles. From the definition it
can be seen that the ranked sizes of the cycles form a
coagulation-fragmentation process (see Aldous (1999) and Pitman
(2002),(2003)) in which components of size $x$ and $y$ merge at
rate $K_n(x,y) = 2xy/n^2$ and components of size $x$ split at rate
$F_n(x) = x(x-1)/n^2$ and are broken at a uniformly chosen random
point. Diaconis, Mayer-Wolf, Zeitouni, and Zerner (2003) have
recently considered the corresponding Markov chain on partitions
of the unit interval and shown that the Poisson-Dirichlet
distribution is the unique invariant measure.

\medskip To study the evolution of the cycles in the random permutation, we
construct a random graph process. Start with the initial graph on
vertices $\{1,\dots,n\}$ with no edge between the vertices. When a
transposition of $i$ and $j$ occurs in the random walk, draw an
edge between the vertices $i$ and $j$. To take care of the rare
event that a given transposition is chosen several times, we will
allow the possibility of multiple edges, and draw a second edge if
one is already present. It is easy to see that in our continuous
time process, at time $t$ this graph is a realization of the
Erd\H{o}s-Renyi random graph $G(n,p)$, in which edges are
independently present with probability $p=1-\exp(-2t/n^2)$, see
Bollob\'as (1985) or Janson, Luczak, and Ruczinski
(2000)\footnote{The fact that we allow multiple edges makes no
difference. At each point where the distinction with the usual
Erd\H{o}s-Renyi random graph may be relevant, a very simple
calculation shows that the effect of multi-edges can be neglected
(see Janson et al. (1993), where this issue is discussed). To make
the core of our arguments simpler to follow, we will ignore this
distinction from now on.}. It is also easy to see that in order
for two integers to be in the same cycle in the permutation it is
necessary that they are in the same component of the random graph.

To estimate the difference between cycles and components, let
$F_t$ denote the event that a fragmentation occurs at time $t$. It
is clear that
\begin{equation}
\label{distance-frag} D_t = N_t - 2 \sum_{s\leq t} \indic{F_s}
\end{equation}
A fragmentation occurs in the random permutation when a
transposition occurs between two integers in the same cycle, so
tree components in the random graph correspond to unfragmented
cycles in the random walk. (Here and in all that follows, "tree"
has a multi-graph meaning : it is a connected component with no
nontrivial closed circuit.) Unicyclic components (with an equal
number of vertices and edges) correspond to cycles in the
permutation that have experienced exactly one fragmentation, but
we need to know the order in which the edges were added to
determine the resulting cycles. For more complex components, the
relationship between the random graph and the permutation is less
clear. Fortunately, these can be ignored in the proofs of our
results.

\section{Limit Theorems}

We will now describe our results and sketch their proofs. Rigorous
proofs of the results stated in this section can be found in
sections 4, 5 and 6.

\subsection{The subcritical regime}

\begin{theorem} \label{subcrit} Let $0<c<1$. The number of fragmentations
\begin{equation}
Z_c := \sum_{s \leq cn/2} \indic{F_s} \cvlaw \text{\em
Poisson}(\kappa (c)) \label{ctsub}
\end{equation}
where $\kappa (c) = (-\log(1-c)-c)/2$. In fact, the convergence
holds for the process $\{ Z_c : 0 \le c < 1\}$ with the limit
being a Poisson process with compensator $\kappa(c)$.
\end{theorem}

\mn {\bf Remark.} The result for fluctuations is formulated in
terms of fragmentations rather than the distance, since $D_{cn/2} -
cn/2 \approx N_{cn/2} - cn/2 = O(n^{1/2})$. For the embedded
discrete time chain, if $k = \lfloor cn/2 \rfloor$, then
\begin{equation}
(k-D_k)/2 \Rightarrow \text{ Poisson($\kappa (c)$) as
$n\to\infty$} \label{dtsub}
\end{equation}
We divide by 2 since a fragmentation reduces the distance by 1
instead of increasing it by 1. To deduce (\ref{dtsub}) from
(\ref{ctsub}) we note that time $k$ in the discrete walk
corresponds to time $N^{-1}(k) \approx cn/2$ in the continuous
time walk.

\mn {\bf Sketch of the proof.} The process $\{Z_c, 0 \le c <1\}$
is a c\`adl\`ag counting process. Therefore by arguments from
Jacod and Shiryaev (1987), it is enough to show that its
compensator $\kappa^n$ converges to the deterministic limit
$\kappa(c)$. If $f_k(t)$ is the fraction of vertices that belong
to cycles of size $k$, the rate at which fragmentations occur is
just $\sum_k f_k(t)(k-1)/n$. Hence $\kappa^n$ is just the integral
with respect to time of this rate. We first show that the variance converges
to 0 and then, by Chebycheff's inequality, it only remains to show
$E \kappa^n(c) \to \kappa(c)$. But by exchangeability
$E[f_k(t)]=P[|\Cc_1|=k]$ where $|\Cc_1|$ is the size of the
component that contains 1 at time $t$. It is not hard to see that
this quantity at time $bn/2$ converges in distribution to the
total progeny $\tau$ of a Galton-Watson branching process with
offspring distribution Poisson($b$), or $PGW(b)$.  Summing the
geometric series, we see that $E\tau = 1/(1-b)$. Integrating
with respect to~$b$ we get the desired expected value, $\kappa(c)$. \qed

\medskip
To prepare for later developments, it is useful to take a second
combinatorial approach to this result. We begin with Cayley's
result that there are $k^{k-2}$ trees with $k$ labeled vertices.
At time $cn/2$ each edge is present with probability $1 -
\exp(-c/n) \sim c/n$ so the expected number of trees of size $k$
present is \beq \binom{n}{k} k^{k-2} \left( \frac{c}{n}
\right)^{k-1} \left( 1 - \frac{c}{n} \right)^{k(n-k) +
\binom{k}{2} - k + 1} \label{Etrees} \eeq since each of the $k-1$
edges needs to be present and there can be no edges connecting the
$k$ point set to its complement or any other edges connecting the
$k$ points. For fixed $k$ the above is asymptotic to
$$
n \frac{k^{k-2}}{k!} c^{k-1} \left( 1 - \frac{c}{n} \right)^{kn}
$$
The quantity in parentheses at the end converges to $e^{-ck}$ so
we have an asymptotic formula for the number of tree components at
time $cn/2$. As a side result we get the following known result:

\begin{cor} \label{borel}
The probability distribution of the total progeny $T$ of a
\emph{Poisson($c$)} branching process with $c <1$ is given by
$P(T=k)= \frac1c \frac{k^{k-1}}{k!} (c\text{e}^{-c})^k $
\end{cor}

See section 4.1 of Pitman (1999) for another proof of this result.
It was first discovered by Borel (1942) and the distribution of
$T$ is called the Borel distribution. It is a particular case of
the so-called Borel-Tanner distribution, see Devroye (1992) and
Pitman (1998) for further references. In this context it appeared
in the problem of the total number of units served in the first
busy period of a queue with Poisson arrivals and constant service
times. See also Tanner (1961). Of course, this becomes a branching
process if we think of the customers that arrive during a person's
service time as their children.

\subsection{The critical regime}

It is well known in the theory of random graphs that the correct time-scale to describe the critical regime is $(n/2)(1 + \lambda n^{-1/3})$,
$\lambda \in (-\infty,\infty)$. See Aldous (1997) for an
interesting account that relates the growth of large clusters in
the critical random graph to the multiplicative coalescent. At
times $(n/2)(1-n^{-r})$ with $r< 1/3$, we are still in the
subcritical regime, so the arguments in the proof of Theorem
\ref{subcrit}, when done more carefully, are still valid. More
precisely, we can show that if $c_n(r)=1-n^{-r/3}$ for $0 \le r
\le 1$, then the expected number of fragmentations up to time
$c_n(r)n/2$ is again given by $\kappa(c_n(r)) \sim (r/6) \log n$.
Hence define:

\begin{equation}
W_n(r) =\left(\frac6{\log n}\right)^{1/2} \left({\ds \sum_{s \leq
c_n(r)n/2}\indic{F_s} - \frac{r}{6} \log n }\right)
\end{equation}

\begin{theorem}
\label{crit1} As $n\to\infty$, $W_n(\cdot)$ converge weakly, with
respect to the Skorokhod topology on the space of c\`adl\`ag
functions on $[0,1]$, to $\{W(r), 0 \leq r \leq 1\}$, a standard
Brownian Motion on $[0,1]$. Furthermore,
\begin{equation}
\left(\frac6{\log n}\right)^{1/2}\left(\sum_{s \leq n/2}
\indic{F_s} - \frac1{6} \log n \right) \Rightarrow W(1),
\end{equation}
\end{theorem}

\mn {\bf Sketch of the proof.} Intuitively, the first result is an
immediate consequence of the Poisson limit in Theorem
\ref{subcrit} and the normal approximation to the Poisson. To
prove it, we show that $W_n(r)$ is a martingale, whose jumps are
asymptotically zero, and whose quadratic variation process is $r$
thanks to our time-change $c_n(r)=1-n^{-r/3}$. Therefore it
converges to Brownian Motion.

At times $(1-n^{-1/3})n/2 \le t \le n/2$ we are in the critical
range of the random graph. Results of Luczak, Pittel, and Wierman
(1994) and computations with (\ref{Etrees}) imply that the number
of fragmentations in this interval is bounded in expectation and
hence can be ignored. \qed

\mn {\bf Remark.} While Theorem \ref{crit1} is a nice theoretical
result, it does not have much to say about any biological example.
If we think of the human genome and set $n = 3$ billion
nucleotides, Theorem 2 says that after $n/2 = 1.5$ billion
transpositions there have been an average of $(\log n)/6 = 3.63$
fragmentations, with a standard deviation of $1.91$. These numbers
are small so even for $n = 3$ billion, we can't expect a very good
approximation to the normal distribution. In the example that we
simulated $n=100$ and $(\log n)/6 = 0.767$ versus an observed
average number of fragmentations = 0.662 (which translates into a
value of 1.224 in Figure 2). While our estimation of the mean is
not very accurate, Figure 3 shows that the distribution of the
number of fragmentations is almost Poisson.

\subsection{The supercritical regime}

This is the most interesting case, and also the hardest one. We
start by establishing a law of large numbers. For all $c >0$
define
$$
\beta_k(c)=\frac1c \frac{k^{k-1}}{k!}(c\text{e}^{-c})^k
$$
so that for $c<1$ it coincides with the Borel distribution of
Corollary \ref{borel}. When $c >1$,
$$
\lim_{n \to \infty} P(|\Cc_1|=k)=\beta_k(c)
$$
still holds but the $\beta_k(c)$'s no longer sum up to 1 because
there is a probability $\beta_{\infty}(c) = 1 - \sum_{k \ge 1}
\beta_k(c)>0$ that $\Cc_1$ is the giant component.

Let us denote by $\Upsilon(c)$ a random variable that takes the
value $1/k$ with probability $\beta_k(c)$ when $1 \le k< \infty$
and the value 0 with probability $\beta_{\infty} (c)$. The
motivation for this definition is that $1/|\Cc_k|$ has the same
distribution as $\Upsilon(c)$ and $\sum_{k=1}^n 1/|\Cc_k|$ gives
the number of components in the random graph.

\begin{theorem}
Let $c>0$ be a fixed positive number. Then the number of cycles in
the random permutation at time $cn/2$, $|\sigma_{cn/2}| = g(c)n +
\omega(\sqrt{n})$, where
\begin{equation}
g(c) := E\Upsilon(c)= \sum_{k=1}^{\infty} \frac1c
\frac{k^{k-2}}{k!} (c\text{e}^{-c})^k
\end{equation}
and the error term $\omega(\sqrt{n})/a_n \sqrt{n} \to 0$ in probability
if $a_n \to \infty$.
\end{theorem}

Note that the theorem is valid for all regimes and implies that
the distance is given by $D_{cn/2} = u(c)n + \omega(\sqrt{n})$
where $u(c)=1-g(c)$. Although it is not obvious from the formula,
$u(c)=c/2$ for $c<1$ and $u(c) < c/2$ when $c > 1$. Using
Stirling's formula, $k! \sim k^k e^{-k} \sqrt{2\pi k}$, it is easy
to check that $g'$ exists for all $c$ and is continuous, but
$g''(1)$ does not exist. In words, there is phase transition in
the behavior of the distance of the random walk to the identity at
time $n/2$  from linear to sublinear.

\begin{proof}
In the supercritical regime the dynamics of the large components
is quite complicated, but there can never be more than $\sqrt{n}$
components of size $\sqrt{n}$ or larger. The expected number of
fragmentations that produce clusters of size smaller than $\sqrt{n}$ by
time $cn/2$ is at most $n^{-1/2} \cdot cn/2$. From this and
Chebyshev's inequality we see that up to a term $\omega(n^{1/2})$,
$|\sigma_{cn/2}|$ is the number of components of the random graph,
and the result follows Theorem 12 in Chapter V of Bollob\'as
(1985).
\end{proof}

\begin{theorem} \label{supercrit-CLT} Let $c>1$. As $n\to\infty$,
\begin{equation}
\label{clt} \frac{D_{cn/2}-u(c)n}{n^{1/2}} \Rightarrow
\mathcal{N}(0,\sigma^2)
\end{equation}
where $\sigma = \rho[1+\rho(c/2-1)]$, and $\rho =1-\theta(c)$ is
the extinction probability of a supercritical $PGW(c)$.
\end{theorem}

\mn {\bf Remark.} Note that the constant $\sigma$ is different
from the one given in Berestycki and Durrett (2003). We were
correct in claiming that the central limit theorem in Theorem
\ref{supercrit-CLT} is the same as the one for the number of
components of the random graph, but we naively thought that the
terms in $\sum_{k=1}^n 1/|\Cc_k|$ were sufficiently independent so
that $\sigma^2 = \var(\Upsilon(c))$.

\mn {\bf Sketch of Proof.} By Pittel's (1990) central limit
theorem for the number of components of a random graph, it suffices
to prove that
the number of extra components due to fragmentation at time $cn/2$
is $o(\sqrt n)$ (see his Corollary 1 and note that $T/c=\rho$).
 Our first step is to increase the cutoff for
large cycles to $n^a$ where $a>1/2$, so that the number of large
cycles is at most $n^{1-a} = o(n^{1/2})$. The number of fragmentations
that produce ``small'' cycles is now $n^{-(1-a)} \cdot cn/2 =
O(n^a)$ and cannot be ignored, so we need to use the fact that
fragmented cycles are reabsorbed by the large components. If the
fraction of mass in large cycles (``upstairs'') at time $tn$ is
$\lambda_t$ then new fragments of size $k$ are produced at rate
$\le 2 \lambda_t$ and each fragment of size $k$ is reabsorbed at
rate $2k \lambda_t$. After time change this is bounded by an
$M/M/\infty$ queue in which the expected number of customers in
equilibrium is $1/k$. Using this, we can show that with high
probability the number of small fragments at any time is at most
$(\log n)^2$. Of course, the coagulation fragmentation process is
not exactly the queuing system. Customers can split into two,
coalesce with other customers, gain weight (and increase their
fragmentation rate) by eating small components, etc. However,
$(\log n)^2$ is much smaller than $n^{1/2}$ so crude but robust
estimates and patience eventually lead to a proof.  \qed

\subsection{Results for Reversals.}

Theorems 3 and 4 extend easily to the approximate distance for
reversal chain. Recall that the main difference lies in the fact
that, a reversal involving edges from different components in the
breakpoint graph always yields a coagulation, but one involving
two edges in the same component may or may not cause a
fragmentation. The proofs of Theorems 3 and 4 for transpositions are
based on showing that fragmentations can be ignored, so this
difference is unimportant and these results extend to reversals. As
Figure 2 shows, this is not true for the more precise results in
Theorems \ref{subcrit} and \ref{crit1}. For example, the
underlying data shows that up to $c=1$, an average of 23\% of the
reversals have caused no change in the distance. Since inversions
that affect an edge are much more frequent than those that involve
it, it seems reasonable to guess that in the limit as $n\to\infty$
the relative orientations of the black edges in a component of the
breakpoint graph are independent. This would imply that the
Poisson process of fragmentations in the reversal case is a
1/2-thinning of the one for transpositions, and Theorem
\ref{crit1} would hold with 6 replaced by 12.

\subsection{Emergence of a giant cycle?}

Since cycles in the random permutation are smaller than components
of the random graph, it follows that if $c<1$ then the largest
cycle at time $cn/2$ has fewer than $\alpha(c)^{-1} \log n$ vertices, where
$\alpha(c) = (c-1-\log c)$. (See Theorem 10 in Chapter V of
Bollob\'as (1985) or Lemma \ref{tailbound} below.)

\medskip For $c>1$, the largest component of the random graph is, as is
well known, ``giant, '' meaning that it is of order $n$. In fact
it is asymptotic to $\theta(c)n$ where $\theta(c)$ is the survival
probability of a supercritical Poisson Galton-Watson with mean
$c$. It is a natural question to ask whether the largest cycle of
the random permutation is also giant in the supercritical regime.

\mn {\bf Conjecture.} {\it Let $L_1(t)$ be the size of the largest
cycle at time $t$. If $c>1$ then
$$
\frac{L_1(cn/2)}{\theta(c)n} \cvlaw V
$$
where $V$ is a random variable with $0<V\le 1\ a.s$.}

\mn This problem is quite different from our original one. However
our techniques enable us to prove a partial result in this
direction as a corollary of the proof of Theorem
\ref{supercrit-CLT}.

\begin{theorem}
\label{upstairs-mass} For any $c>1$, at time $cn/2$ there are at
least $\theta(c) n -o(n)$ vertices located on large cycles (i.e.,
of size  greater than or equal to $n^a$, for any $a<2/3$).
\end{theorem}

\mn David Aldous (private communication) conjectures that the
relative sizes of the pieces of the giant cycle are in equilibrium
at all times in the supercritical regime, i.e., have the
Poisson-Dirichlet $PD(0,1)$ distribution, which gives the limiting
behavior of the ordered sizes of cycles in a uniform random
permutation. According to this conjecture, $V$ would be
distributed as the first coordinate of a $PD(0,1)$ random
variable. One way to approach this conjecture would be to
generalize Aldous (1997) to show that the large cycles in the
critical regime converge to a coagulation-fragmentation process
and to study the growth of clusters in that process.

\medskip
 Alternatively, one could look at the size of the cycle
containing 1, $K_1(t)$, and try to show that
$$
\frac{K_1(cn/2)}{\theta(c)n} \cvlaw U
$$
where $U$ has a point mass of size $1-\theta(c)$ at 0 and is
otherwise uniform on $(0,\theta(c))$. Figure 4 shows the average
growth of $K_1(cn/2)/n$ in 10,000 simulations of $n=100$,
$n=1000$, and compares the results to $EU=\theta(c)^2/2$. Although
this considers only one aspect of the distribution of large
cycles, it agrees well with Aldous' conjecture.

\mn Figure 5 shows a histogram of the result of 100,000
simulations of $K_1(100)$ when $n=100$. As the graph shows, the
spike in the frequency of clusters of size 4 or smaller is what
one would predict from the random graph cluster size distribution.
The remainder of the distribution is roughly uniform except for
rounding at the upper end. The latter is to be expected if Aldous'
conjecture is correct, since the size of the giant component
satisfies the central limit theorem.

\medskip As we were finishing this paper, we learned that Oded
Schramm (private communication) has proved David Aldous'
conjecture.

\mn \textbf{Remark.} The problem of the emergence of a giant cycle
is closely related to Angel's (2003) work on the existence of
infinite orbits for the \emph{random stirring process}, which is
the random transposition random walk on an infinite graph such as
$\mathbb{Z}^d$ or a tree, rather than the complete graph on
$\{1,\ldots, n\}$ considered in this work. To explain the
connection, suppose that we construct our process using a Poisson
process with rate $2/n^2$ for each $i\neq j$, and at these times
draw an edge between $i$ and $j$ to indicate that $i$ and $j$ are
to be transposed. To compute the cycles in the permutation at time
$cn/2$, we repeat the first $[0,cn/2]$ units of time periodically
and then observe the sites that a walker starting at $i$ visits at
times $kcn/2$, for $k=1,2, \ldots$. Angel (2003) calls this
construction the \emph{cyclic time random walk}. Its relevance to
his work is that the cyclic time random walk is transient if, and
only if, the cycles are infinite.

\section{The subcritical regime}

Let us introduce some notations for the different probability laws
involved. For each $n$, we have the coagulation-fragmentation
process, and the Erd\H{o}s-Renyi random graph model. To emphasize
when computations are being done for the random graph we will use
$Q_p$, for the random graph with Bernoulli percolation parameter
$p$, and $Q$ for the law of the evolving random graph that at time
$s$ has $p_s = 1 - \exp(-2s/n^2)$. When $s=cn/2$ this probability
is $p(c,n) = 1 - \exp(-c/n) \le c/n$. To simplify notation we will
use $Q X$ to denote the expected value of $X$ with respect to the
probability $Q$.

\subsection{Preliminary results : comparison with a branching process}

Our first result provides a useful upper bound.

\begin{lemma}
The cluster size $|\Cc_1|$ in $Q_{c/n}$ is dominated by $Z$, the
total progeny of a branching process in which each individual has
a Binomial$(n-1,c/n)$ number of children, i.e., we can construct
these random variables on the same probability space so that
$|{\Cc}_1| \le Z$ a.s. It follows from this that if $c<1$ then
$Q_{c/n}|{\Cc}_1| \le 1/(1-c)$. \label{dominate}
\end{lemma}

\begin{proof}
Intuitively, this holds since a vertex in generation $k$ may have
children among all of the $n$ vertices of the graph except those
of the first $k$ generations. To begin to prove this formally, let
$\xi_{i,j}$, $1 \le i,j\le n$ be independent random variables,
taking values 1 with probability $c/n$ and 0 with probability
$1-c/n$. To start the random graph let $Y_0=\{1\}$ and let $Y_1 =
\{ j \not\in Y_0: \xi_{1,j} = 1 \}$. To start the branching
process let $Z_0=1$, $Z_1 = |Y_1|$, and let $\phi_1: Y_1 \to \{ 1,
2, \ldots Z_1\}$ be 1-1 and onto.

If the first $k$ stages of the construction have been done and we
have $Y_{k} \neq \emptyset$ and a $\phi_k: Y_k \to \{1, \ldots Z_k
\}$ that is 1-1 (but not onto in general), then let
$$
Y_{k+1} = \cup_{i \in Y_k} \{ j \not\in \cup_{\ell=0}^k Y_\ell:
\xi_{i,j} = 1 \}
$$
We let individual $\phi_k(i)$ in the $k$th generation of the
branching process have $|\{ j \neq i: \xi_{i,j} = 1 \}|$ children.
The individuals in the branching process that are not in
$\phi_k(Y_k)$ have a number of children given by independent
binomials. It should be clear from the construction that can again
define $\phi_{k+1}: Y_{k+1} \to \{1, \ldots Z_{k+1} \}$ to be 1-1,
and the comparison follows by induction. The inequality follows by
computing $EZ$ (for instance by summing a geometric series).
\end{proof}

The next result shows that the bound in Lemma \ref{dominate} is
exact in the limit. Let $\{Z_k\}_{k=0}^{\infty}$ be a Poisson
Galton-Watson process with offspring mean $c$ and let
$Z=\sum_{k=0}^{\infty} Z_k$ be its total progeny.

\begin{lemma} \label{PGW-limit}
Let $\Cc_1$ be the cluster that contains vertex 1. If $0 \le c <
1$ then as $n\to\infty$
$$
Q_{p(c,n)} (|\Cc_1| = k ) \to P( Z = k )
$$
\end{lemma}

\begin{proof}
The number of children of vertex 1, $Z^1_n=|Y_1|$ has distribution
Binomial$(n-1,p(c,n))$, which converges to a Poisson(c) limit. Let
$k \geq 1$ and let $(n_1,...,n_{k+1}) \in \Nn^{k+1}$. If we let
$Z_j^n = |Y_j|$ then
$$
Q_{p(c,n)}(Z_{k+1}^n=n_{k+1}|Z_1^n=n_1,...,Z_k^n=n_k) = P
\left(\sum_{i=1}^{n_k} B_i^n=n_{k+1}\right)
$$
where $B_i^n$ are i.i.d.~Binomial$(n-s,p(c,n))$ random variables,
and $s=\sum_{i=0}^k n_k$ with $n_0=1$. From this it follows easily
that the convergence of finite-dimensional distributions of
$\{Z_j^n\}_{j \geq 1}$ to those of $PGW(c)$. Markov's inequality
and the domination result in Lemma \ref{dominate} imply that
$$
Q_{p(c,n)}\left(\sum_{k=K}^\infty Z^n_k > 0 \right) \le
Q_{p(c,n)}\left(\sum_{k=K}^\infty Z^n_k\right) \le c^K/(1-c)
$$
and the desired conclusion follows.
\end{proof}

Our next ingredient is

\begin{lemma}
\label{tailbound} $Q_{c/n}(|\Cc_1| \ge y ) \le c^{-1} \exp( - (c-1
- \ln c) y )$.
\end{lemma}

\begin{proof}
In view of Lemma \ref{dominate}, it suffices to prove the result
for $Z$, rather than $|\mathcal{C}_1|$. To do this, let \beqax
\phi_n(\theta) & = & e^{-\theta} \sum_{m=0}^{n-1} {n-1 \choose m}
\left( \frac{c}{n} \right)^m
\left( 1 - \frac{c}{n} \right)^{n-1-m} e^{\theta m} \\
& = & e^{-\theta} \left( 1- \frac{c}{n} + \frac{c}{n} e^\theta
\right)^{n-1} \eeqax be the moment generating function of the
distribution of the number offspring minus 1. Let $S_m$ be a
random walk that takes steps with this distribution and $S_0=1$,
so that $S_m$ explores the Galton-Watson tree. Then $\tau = \inf\{
m: S_m = 0\}$ has the same distribution as $Z$. Let $R_m =
\exp(\theta S_m)/\phi_n(\theta)^m$. $R_m$ is a nonnegative
martingale. Stopping at time $\tau $ we have $e^{\theta} \ge
E(\phi_n(\theta)^{-\tau})$. If $\phi_n(\theta)<1$ it follows that
$$
P(\tau \ge y ) \phi_n(\theta)^{-y} \le E[\phi_n(\theta)^{-\tau}]
\le e^{\theta}
$$
Using $\phi_n(\theta) \le e^{-\theta} \exp(c(e^\theta - 1))$ now
we have
$$
P(\tau \ge y ) \le e^{\theta} \left(e^{-\theta}\exp(c(e^\theta -
1))\right)^y
$$
To optimize the bound we want to minimize $c(e^\theta - 1) -
\theta$. Differentiating this means that we want $ce^\theta - 1 =
0$ or $\theta = -\log(c)$.  Plugging this and recalling that
$\tau$ and $Z$ have the same distribution we have
$$
P(Z \ge y ) \le \frac1c \exp( - (c-1 - \ln c) y )
$$
It follows that
$$
Q_{c/n}(|\Cc_1| \ge y ) \le \frac1c \exp( - (c-1 - \ln c) y )
$$
which completes the proof of Lemma \ref{tailbound}.
\end{proof}

Now recall that for $c<1$, $Z_c= \sum_{s\le cn/2}\indic{F_s}$ is
the number of fragmentations up to time $cn/2$.

\begin{lemma} \label{PGW-limitb}
Let $f_k(s)$ be the empirical fraction of vertices in cycles of
size $k$ at time $s$. If $0 \le c < 1$ then $E f_k(cn/2) \to P( Z
= k )$ and $E Z_c \to \kappa(c)$, where $\kappa(c)$ was defined in
Theorem \ref{subcrit}.
\end{lemma}

\begin{proof}
The cycle sizes at time $s$ in the coagulation-fragmentation
process are dominated by the cluster sizes in the random graph
model with $p_s = 1 - \exp(-2s/n^2) \le 2s/n^2$. Therefore,
$$
E Z_c \le \int_0^{cn/2} Q f_k(s) \frac{k-1}{n} \, ds \le
\int_0^{cn/2} Q_{2s/n^2} |\Cc_1| - \frac{1}{n} \, ds
$$
Using Lemma \ref{dominate} $Q_{2s/n^2} |\Cc_1| \le 1/(1- (2s/n))$.
Changing variables $un/2 = s$ we have \beq E Z_c \le -
\frac{1}{2}(\log(1-c) + c) = \kappa(c)  \label{bdfrag} \eeq Since
unfragmented cycles are the same as tree components in the random
graph, the first convergence result follows from Lemma
\ref{PGW-limit}. The second one follows from Fatou's lemma and
(\ref{bdfrag}).
\end{proof}

The final preparatory step is:

\begin{lemma} \label{twofrags}
If $c<1$ the expected number of fragmentations that occur to
cycles that have already been fragmented is $\le  K_c (\log
n)^2/n$, and $K_c = 9c \kappa(c) \alpha(c)^{-2}$. (Recall
$\alpha(c)= (c-1-\log c)$).
\end{lemma}

\begin{proof}
The expected number of such fragmentations is at most:

\beqax
\le &\ds E\int_0^{cn/2} \frac{\# \text{vertices in fragments}}n\frac{L_1(bn/2)}n dt\\
 & \ds \frac{n}2 \int_0^c E Z_b \left(\frac{L_1(bn/2)}n\right)^2 db
\eeqax where $L_1(t)$ is the size of the largest component at time
$t$. In the event that $L_1(cn/2) \le 3\alpha(c)^{-1} \log n$, the
above is at most
$$(n/2)(3\alpha(c)^{-1}\log n/n)^2\int_0^c \kappa(b)db \le \frac12K_c\frac{(\log n)^2}n$$
On the other hand by Lemma \ref{tailbound} the complement of this
event has probability at most $n^{-2}$, and there can never be
more than $cn/2$ such fragmentations, so Lemma \ref{twofrags} is
proved.
\end{proof}

\subsection{Proof of Theorem 1}
We are now ready to prove Theorem \ref{subcrit}. Let $\bar Z^n_c =
\sum_{s\le cn/2} \indic{\bar F_s}$, $0 \le c < 1$ be the counting
process of fragmentations that occur to cycles which (a) have not
been fragmented previously and (b) have size $\le n^{0.7}$. The
second condition is irrelevant in this section, but imposing it
now will help in the next one. Unfragmented cycles correspond to
trees in the random graph so the compensator of $\bar Z^n_c$ is
\begin{equation}
\bar\kappa^n(c) = \int_0^{cn/2} \bar\psi^n_s \, ds
\end{equation}
where $\bar\psi^n_s = \sum_{k=1}^{n^{0.7}} \bar f_k(s) (k-1)/n$
and $\bar f_k(s)$ is the fraction of vertices that belong to tree
components of size $k$. As noted in the sketch of the proof, it is
enough to show that for each fixed $c$, $\kappa^n(c)$ converges in
probability to $\kappa(c)$, or, by Lemma \ref{twofrags}, that
$\bar{\kappa}^n(c)$ converges to $\kappa(c)$ in probability.
Lemmas \ref{tailbound} and \ref{PGW-limitb}  imply that
$E[\int_0^{cn/2}\bar{\psi}_s^n ds] \tendsto \kappa(c)$. It remains
to show that $\var \int_0^{cn/2}\bar\psi_s^n ds \tendsto 0$. Our
first step will be to prove :
\begin{equation}
\label{variancebound}
 \mbox{\rm
var}\,(\bar\psi^n_s) \leq \frac{K}{n^3} Q_{p(c,n)}[|{\Cc}_1|^3]
\end{equation}
for all time $s \leq cn/2$, where $K$ is a constant that depends
only on $c$.

\medskip To see this, first observe that in terms of cluster sizes
$$
\bar\psi^n_s = \frac{1}{n^2}\sum_{i=1}^n (|\Cc_i|-1)I_i
$$
where $I_i$ is the indicator of the event that $\Cc_i$ is a tree.
Let $d_i = (|\Cc_i|-1)I_i$. \beq \var \frac1{n^2}(d_1 + \cdots +
d_n ) = \frac1{n^4} \left( n \var(d_1) + n(n-1) \cov(d_1,d_2)
\right) \label{varf} \eeq Monotonicity and Lemma \ref{tailbound}
imply, \beq \var (d_1) \leq Q_{p(c,n)}[|\Cc_1|^2] \leq K
\label{varb} \eeq It remains to bound $\cov (d_1,d_2)$. If we let
$$
\pi_{n,i} = i^{i-2} \left(p_s\right)^{i-1} \left(1 -
p_s\right)^{i(n-i) + {i\choose 2} - (i-1)}
$$
where $p_s = 1 - \exp(-2s/n^2)$ then by the reasoning for
(\ref{Etrees}) we have \beqax &&Q_{p_s}[ {\Cc}_1 \cap {\Cc}_2 =
\emptyset, |{\Cc}_1| = j, |{\Cc}_2| = k,\ \Cc_1 \text{ and } \Cc_2
\text{ are trees}] =
{n-2 \choose j-1} \pi_{n,j} {n-j-1\choose k-1} \pi_{n-j,k} \\
&&Q_{p_s}[ {\Cc}_1 = {\Cc}_2, |{\Cc}_1| = k, \ \Cc_1 \text{ is a
tree}] = {n-2 \choose k-2} \pi_{n,k} \eeqax From this it follows
that $\cov(d_1,d_2)$
 \beqax & = &\ds \sum_{j,k} \left[ {n-2
\choose j-1} {n-j-1\choose k-1}(1-p_s)^{-j}
- {n-1 \choose j-1} {n-1\choose k-1} \right] (j-1)(k-1) \pi_{n,j} \pi_{n,k} \\
 & & \ds + \sum_k {n-2 \choose k-2} (k-1)^2 \pi_{n,k}
\eeqax
For the first term in the right-hand side,
\beqax
&& \ds\left[ {n-2 \choose j-1} {n-j-1\choose
k-1}(1-p_s)^{-j}-{n-1 \choose j-1} {n-1\choose k-1} \right]   \\
&& \le \ds \frac{(n-2)! e^c}{(j-1)!(k-1)!(n-j-k)!}
- \frac{(n-1)!}{(j-1)!(n-j)!} \frac{(n-1)!}{(k-1)!(n-k)!}  \\
 && \le 0
\eeqax

\mn since $(n-2)! e^c \le (n-1)!$ for large $n$ and
$(n-j)!/(n-j-k)! \le (n-1)!/(n-1-k)!$.

\mn For the second term,
$$\sum_k {n-2 \choose k-2} (k-1)^2
\pi_{n,k}
  \leq \frac{1}{n-1} \sum_k k^3 {n-1 \choose k-1}
\pi_{n,k}
  \leq \frac{1}{n-1} Q_{p(c,n)}[|{\Cc}_1|^3] $$
 Combining this
with (\ref{varf}) and (\ref{varb}) gives (\ref{variancebound}).

\mn Hence by the Cauchy-Schwarz inequality we get: \beqa \var
\left( \int_0^t \bar\psi^n_s ds \right) &=& Q \left[\left(\int_0^t
(\bar\psi_s^n - Q[\bar\psi_s^n])\, ds \right)^2\right] \leq t
\int_0^t \var(\bar\psi_s^n) ds
\label{CS}\\
&\leq & \ds \frac{cn}2 \int_0^{cn/2}\frac{K}{n^3} ds =
\frac{c^2K}{4n} \tendsto 0 \nonumber \eeqa where we have used both
(\ref{variancebound}) and Lemma 3. This completes the proof of
Theorem \ref{subcrit}.

\section{The critical regime}

The first step in the proof of Theorem \ref{crit1} is to argue
that fragmentations of previously fragmented cycles can be
ignored. The number of such fragmentations is smaller than the
total number of cycles in multicyclic components (i.e., components
with at least 2 cycles) in the random graph. Theorem 1 and
Corollary 3 in Luczak, Pittel, and Weirman (1994) imply that the
total number of cycles in multicyclic components in the critical
regime is bounded in probability.\footnote{This result could also
be derived from the Folk Theorem 1 in Aldous (1997) which gives
the limit for the joint distribution of the component sizes and
the number of cycles they contain. See the discussion page 850 of
his paper.} In particular, divided by $(\log n)^{1/2}$ it
converges to 0 in probability. As a result, by the converging
together lemma (see e.g., Durrett (1996), Chap.2, Ex. 2.10), it
suffices to prove the central limit theorem for the number of
fragmentations on tree components.

As in the previous section, we will in addition restrict our
attention to fragmentations of tree components of size at most
$n^{0.7}$, and continue to use the notation introduced there.
(Indeed, classical results from the theory of random graphs, or
Aldous (1997), show that asymptotically almost surely all clusters
are smaller than $n^{0.7}$).

\mn Let $\bar{W}_n(r):= (6/ \log n)^{1/2} (\bar{Z}^n(r) -
\bar{\kappa}^n(r))$. By standard methodology in the theory of
stochastic processes (see Jacod and Shiryaev (1987) or Revuz and
Yor(1999) for instance), to prove convergence of
$\bar{W}_n(\cdot)$ to Brownian Motion, the two things we need to
check are: (i) $E[{\sup}_{0 \leq r \leq 1}
|\bar{W}_n(r)-\bar{W}_n(r^-)|] \tendsto 0$ and (ii) The quadratic
variation of $W_n$, i.e. the increasing process associated with
$\bar{W}_n(\cdot)^2$, must converge to $r$ at time $r$. (i) is
obvious because $\bar{Z}^n$ is a counting process, and (ii) turns
into $E(6\bar{Z}^n(r)/\log n) \tendsto r$ and
$\var(6\bar{Z}^n(r)/\log n) \tendsto 0$. These two steps are dealt
with respectively in lemmas \ref{critexp} and \ref{critvar}.

\mn But first, we need a technical lemma that will be useful on
several occasions (e.g., for computing precise asymptotics of the
number of trees of a given size).

\begin{lemma}
\label{Stirling} If $k\to\infty$ and $k=o(n^{3/4})$ then \beqax
\gamma_{n,k}(c) & \equiv & \binom{n}{k} k^{k-2} \left( \frac{c}{n}
\right)^{k-1}
\left( 1 - \frac{c}{n} \right)^{kn - k^2/2 - 3k/2 + 1} \\
& \sim &  \frac{n k^{-5/2}}{c\sqrt{2\pi}} \exp\left( -\alpha(c)k +
(c-1)\frac{k^2}{2n} - \frac{k^3}{3n^2}\right) \equiv
\lambda_{n,k}(c) \eeqax where $\alpha(c) = c-1-\log(c)$. There is
a constant $K$ so that if $1\le k \le n^{0.7}$ and $c \le 1$ then
$\gamma_{n,k}(c) \le K \lambda_{n,k}(c)$.
\end{lemma}

\begin{proof}
Stirling's formula implies $k! \sim k^k e^{-k} \sqrt{2\pi k}$.
Using this we have that
$$
\gamma_{n,k}(c) \sim \frac{n k^{-5/2}}{c\sqrt{2\pi}} \left[
\prod_{j=1}^{k-1} \left( 1 - \frac{j}{n} \right) \right]  e^k c^k
\left( 1 - \frac{c}{n} \right)^{kn - k^2/2 - 3k/2 + 1}
$$
Using the expansion $\log(1-x) = -x - x^2/2 - x^3/3 - \ldots$ we
see that if $k=o(n)$ then
$$
\left( 1 - \frac{c}{n} \right)^{kn - k^2/2 - 3k/2 + 1} \sim \exp(
-ck + k^2/2n )
$$
while if $k=o(n^{3/4})$ we have \beqax \prod_{j=1}^{k-1} \left( 1
- \frac{j}{n} \right) & = & \exp\left( - \frac{1}{n}
\sum_{j=1}^{k-1} j - \frac{1}{n^2} \sum_{j=1}^{k-1} j^2
+ O\left( \frac{j^4}{n^3}\right) \right) \\
& \sim & \exp\left( - \frac{k(k-1)}{2n} -
\frac{k(k-1)(2k-1)}{6n^2} \right)
 \sim \exp\left( - \frac{k^2}{2n} - \frac{k^3}{3n^2} \right)
\eeqax Combining the last three formulas gives the asymptotic
formula. To prove the bound we note that Stirling's formula
implies $k! \ge \delta k^k e^{-k} \sqrt{2\pi k}$ for some
$\delta>0$. Using the bounds $\log(1-x) \le -x$ and $\log(1-x) \le
-x-x^2/2$ in the last two calculations gives the upper bound.
\end{proof}

\begin{lemma}
\label{critexp}
$$
E\left[ \frac{6}{\log n} \int_0^{ c_n(r)n/2} \bar{\psi}_s^n \ ds
\right] \tendsto r
$$
\end{lemma}

\begin{proof}
The upper bound follows from (\ref{bdfrag}) which holds for all
$c<1$. In the other direction, changing variables $s=c_n(v)n/2$
where $c_n(v)=1-n^{-v/3}$ and noting $c_n'(v) = (1/3)(\log
n)n^{-v/3}$ gives \beqa {E\left[\frac{6}{\log n} \bar{Z}^n(r)
\right]} & = & n  \int_0^r Q[\bar{\psi}^n_{c_n(v)n/2}] n^{-v/3} \,
dv
\label{chvar}\\
& = &  \int_0^r \sum_{k=1}^{n^{0.7}} \frac{k-1}n
Q_{p(c_n(v),n)}[kT_k]n^{-v/3} \, dv \nonumber \eeqa where $T_k$ is
the number of tree components of size $k$, and $p(c,n) = 1 -
\exp(-c/n)$.

We can take the limit of the last expression by using formula
(\ref{Etrees}), combined with Lemma \ref{Stirling}. Indeed formula
(\ref{Etrees}) shows that $ET_k =\gamma_{n,k}(c)$, and $k\le
n^{0.7} = o(n^{3/4})$, so that the use of Lemma \ref{Stirling} is
justified. Hence \beqa {E\left[\frac{6}{\log n} \bar{Z}^n(r)
\right]} = \int_0^r \sum_{k=1}^{n^{0.7}} \frac{k(k-1)}n
\gamma_{n,k}(c_n(v)) n^{-v/3} \, dv \nonumber \eeqa

Setting $c=1-b$ with $b=n^{-v/3}\to 0$ and using Taylor's theorem
$$
-(c-1-\log(c))k - b\frac{k^2}{n} = -\frac{b^2}{2}k -
b\frac{k^2}{n} +o(b^2k)\label{asyalpha}
$$
The first term becomes significantly negative when $k \approx
1/b^2=n^{2v/3}$, the second when $k \approx \sqrt{n/b}=n^{1+v/3}$.
When $v\le r<1$ the first threshold is smaller and the second term
can be ignored. Thus Lemma \ref{Stirling} and the last observation
imply that if $v<1$ \beq \sum_{k=1}^{n^{0.7}} \frac{k^2}n
Q_{p(c_n(v),n)}[T_k] \sim \frac{1}{\sqrt{2\pi}}
\sum_{k=1}^{n^{0.7}} k^{-1/2} \exp(-n^{-2v/3}k/2) \label{asyQ}
\eeq Here we have used the asymptotic formula of Lemma
\ref{Stirling} for all $k$. However, the next computation will
show that the sum grows like $n^{v/3}$ so the contributions from
small $k$ can be ignored.

If we view the sum in (\ref{asyQ}) as a Riemann sum with spacing
$n^{-2v/3}$, we can rewrite it as
$$
n^{v/3} \sum_{k=1}^{n^{0.7}} n^{-2v/3} (n^{-2v/3}k)^{-1/2}
\exp(-n^{-2v/3}k/2)
$$
From this it follows that
$$
n^{-v/3} \sum_{k=1}^{n^{0.7}} \frac{k(k-1)}n Q_{p(c_n(v),n)}[T_k]
\to \int_0^\infty \frac{x^{-1/2}}{\sqrt{2\pi}} e^{-x/2} \, dx
$$
Changing variables $x=y^2$, $dx = 2y\, dy$ the integral becomes
$(2\pi)^{-1/2} \int_0^\infty 2 e^{-y^2/2} \, dy = 1$. Therefore,
by Fatou's lemma:
$$
{\liminf}_{n\to\infty} E\left[ \frac{6}{\log n}
\bar{Z}^n(r)\right]
 \ge \int_0^r n^{-v/3}\cdot n^{v/3}\ dv = r
$$
\end{proof}

We turn now to the analysis of the variance.

\begin{lemma}
\label{critvar} $\hbox{\rm var}\left(\frac6{\log n}
\bar{\kappa}_n(r)\right) \tendsto 0$
\end{lemma}

\begin{proof}
Changing variables as in (\ref{chvar}) and using Cauchy-Schwarz
inequality as in (\ref{CS}), \beqax \var \left( \frac{6}{\log n}
\int_0^{ c_n(r)n/2} \bar{\psi}_s^n \ ds \right)
& = &  \var \left( n \int_0^r \bar{\psi}^n_{c_n(v)n/2} \, n^{-2v/3} \, dv \right) \\
& \leq & n^2 \int_0^r n^{-2v/3} \var (\bar{\psi}^n_{c_n(v)n/2}) \, dv \\
& \leq & \frac{2}{n} \int_0^r n^{-2v/3}
Q_{p(c_n(v),n)}[|\Cc_1|^3I_1] \,  dv \eeqax Reasoning as in
(\ref{asyQ}) but using the bound in Lemma \ref{Stirling}
$$
\sum_{k=1}^n k^3 Q_{p(c_n(v),n)}[kT_k] \le K \sum_{k=1}^n k^{3/2}
\exp(-n^{-2v/3}k/2)
$$
To check the right-hand side note that the power of $k$ has
increased by 2, from the previous calculation. If we view the last
sum as a Riemann sum with spacing $n^{-2v/3}$, we can rewrite it
as
$$
n^{5v/3} \sum_{k=1}^n n^{-2v/3} (n^{-2v/3}k)^{3/2}
\exp(-n^{-2v/3}k/2)
$$
Now $x^{3/2}e^{-x/2}$ has derivative $((3/2)x^{1/2} -
(1/2)x^{3/2})e^{-x/2}$ so it is increasing on $[0,3]$ and then
decreasing on $[3,\infty)$. Thus if we discard the term with the
largest $k$ so that $n^{-2v/3}k \le 3$ we have a lower bound on
the integral.
$$
n^{-2v/3} \sum_{k=1}^n k^3 Q_{p(c_n(v),n)}[kT_k]  \le n^v
\frac{1}{\sqrt{2\pi}} \int_0^\infty x^{3/2} e^{-x/2} \, dx  +
n^{v/3}3^{3/2}e^{-3/2}
$$
Using this it follows that
$$
\var \left(\frac6{\log n} \bar{\kappa}^n(r)\right) \leq
\frac{K}{n} \int_0^r n^u \, du
$$
Writing $n^u = \exp(-u \log n))$ and integrating we have that the
right-hand side is $\leq K/(\log n)\to 0$. This concludes the
proof of the first result in Theorem \ref{crit1}.
\end{proof}

The final step is to estimate the number of fragmentations that
occur to tree components of size $\le n^{0.7}$ at times between
$(1-n^{-1/3})n/2$ and $n/2$:
$$
\int_{(1-n^{-1/3})n/2}^{n/2} Q \bar\psi_s^n \, ds
$$
For each $s$ in the interval the integrand is smaller than
$\sum_{k=1}^{n^{0.7}} \frac{k^2}{n} Q_{1/n}T_k$. Using Lemma
\ref{Stirling}, the last quantity is smaller than
$$
\frac{K}{n} \sum_{k=1}^{\infty} k^{-1/2} \exp(-k^3/3n^2)
$$
which we can rewrite as
$$
\frac{K}{n}n^{1/3} \sum_{k=1}^{\infty} n^{-2/3} (k
n^{-2/3})^{-1/2} \exp(-(k n^{-2/3})^3/2)
$$
The above sum is a Riemann sum so it converges to $\int_0^{\infty}
x^{-1/2}e^{-x^3/2}dx $. Therefore, $Q\bar{\psi}^n_s \le K
n^{-2/3}$. Since the duration of the critical regime is
$n^{2/3}/2$, the expected number of fragmentations is bounded and
the proof of Theorem \ref{crit1} is complete.

\section{The supercritical regime}

By Pittel's (1990) central limit theorem for the number of
components of a supercritical random graph, it is enough to show
that, with probability going to 1 as $n \to \infty$, at time
$cn/2$ there are fewer than $o(n^{1/2})$ extra components due to
fragmentation. (This was already indicated in the sketch of the
proof of Theorem \ref{supercrit-CLT}).

\medskip Let $a=0.55$. (In fact the results stated in this section
would also be valid for any $1/2 < a < 2/3$ but making this choice
makes some proofs slightly easier). We call cycles of size $k \ge
n^a$ large. These can be ignored since there cannot be more than
$n^{1-a}=o(n^{1/2})$ such components. We define the amount of mass
``upstairs" by
$$
N_t^{\uparrow} = \sum_{k>n^a} kX_k(t)
$$
where $X_k(t)$ is the number of cycles of size $k$ at time
$n/2+t$. (It is convenient in this section to shift the time so
that $t=0$ corresponds to critical time $n/2$.) If all of the mass
was upstairs, then the expected number of cycles of size less than
$n^a$ produced by fragmentation would be $2n^{a-1} (cn/2) =
O(n^a)$. It is overly pessimistic to think that all of the mass
will be upstairs, but by analogy with the random graph, we expect
(and will eventually prove in Theorem \ref{upstairs-mass}) that at
times $c>1$ a positive fraction of the total mass $n$ will be
there, so this estimate of the number of fragmentations is too
large to ignore.

\medskip To improve this crude estimate, we take advantage of the
fact that fragmented pieces are reabsorbed upstairs. Let
$X^{\downarrow}_k(t)$ be the number of cycles of size $k$ produced
by fragmentation of cycles upstairs. $X^{\downarrow}_k(t)$ can
only increase when a transposition is performed, and only if it is
made of one of the $N_t^{\uparrow}$ vertices upstairs and of one
of the 2 points located $k$ steps away when writing the
corresponding cycle of the current permutation. This gives a rate at
most $2N_t^{\uparrow}/n^2$. As for the death rate, one way to get
rid of a component of size $k$ is by picking one of the $k$
vertices of one of the $X^{\downarrow}_k(t)$ components and one of
the $N_t^{\uparrow}$ vertices upstairs. This happens with rate
$2kX^{\downarrow}_k(t)N_t^{\uparrow}/n^2$. For the moment we are
ignoring the fact that cycles may experience coalescence or
fragmentation while downstairs. We will deal with these
complexities once we have an understanding of the basic birth and
death process of fragments of large clusters.

\subsection{The cluster queuing system}
\label{cqs} It is fortunate that the unknown quantity
$N_t^{\uparrow} \le n$ appears in both rates, so that as along as
$N_t^{\uparrow}>0$ we can remove it by time change. Once this is
done, we have a system of stochastic processes $\xi^k_t$, for
$1\le k \le n^a$ that we call a \emph{cluster queuing system}: let
$\xi^k_t$ be independent birth-and-death chains with birth rate 1
and death rate $k\xi^k_t$, that begin with $\xi_0^k=0$.

\begin{lemma}
\label{bounds-system} With probability $\to 1$ as $n\to\infty$ we
have
$$
\sum_{k=1}^{n^a} \xi^k_t \le (\log n)^2\quad\hbox{and}\quad
\sum_{k=1}^{n^a} k\xi^k_t \le n^a(\log n)^2
$$
for all $t \le c$ ($c>0$).
\end{lemma}

\mn \textbf{Remark.} Although this system of stochastic processes
can be defined without any reference to our random walk problem,
it is useful to bear in mind that the state of this cluster
queuing system at time $t$ describes the number of fragments of
large cycles at time
$$
\frac{n}2 +\int_0^t \frac{n^2}{2N^{\uparrow}_s}ds \ge \frac{n}2(1+t)
$$
since $N^{\uparrow}_s\le n$. Thus the control obtained in the
above lemma for all $t\le c$, will provide useful information for
the random walk between times $n/2$ and $(1+c)n/2$ for any $c>0$.
On our original time-scale, this corresponds exactly to the
supercritical regime, i.e. up to time $cn/2$ for any $c>1$.

\begin{proof}
The second result is a trivial consequence of the first. The key
idea to handle the processes $\xi^k_t$ is to consider strips $2^j
\le k < 2^{j+1}$. Because there are no simultaneous jumps, we can
prove that the queues $\xi^k_t$ at each level $k$ are independent
processes (see e.g. Revuz-Yor (1999), chap. XII, prop. (1.7), for
a proof of this fact in the case of Poisson processes). Therefore,
for each $1 \le j \le \log_2 n^a $, the number of cycles with
sizes in $[2^j, 2^{j+1})$, $\zeta^j_t$, is dominated by a birth
and death chain with birth and death rates respectively $2^j$ and
$2^j\zeta^j_t$. To analyze these processes, we consider the
successive excursions away from 0. Their embedded discrete time
processes $Y_s$ jump from $m$ to $m-1$ with probability $m/(m+1)$
and from $m$ to $m+1$ with probability $1/(m+1)$. Let us try to
find a function $\phi$ such that $\phi(0)=0$, $\phi(1)=1$ and
$\phi(Y_s)$ is a martingale. The latter implies
$$
\frac{1}{m+1}[\phi(m+1)-\phi(m)] =
\frac{m}{m+1}[\phi(m)-\phi(m-1)]
$$
so $\phi(x) = \sum_{k=1}^x (k-1)!$. Since $\phi(1)=1$ and
$\phi(0)=0$, it follows by optional sampling that the maximum
level reached during an excursion of $\zeta^j$, $M$, satisfies
\beq P(M > x) = 1/\phi(x+1) \le 1/x! \label{maxbd} \eeq

To bound the number of excursions for the process in the
$j^{\text{th}}$ strip before time $c$, $N_j(c)$, we note
that jumps from 0 to 1 occur at rate $2^j$ so ignoring the amount
of time it takes to return to 0 from 1, the number of excursions
by time $c$ is bounded by a Poisson random variable with mean
$2^jc \le cn^{a}$. Markov's inequality implies that $P(
N_j(c) > n^2 ) \le cn^{a-2}$ so \beq P\left( \max_{1\le j \le
a\log_2 n} N_j(c) > n^2 \right) \to 0 \label{excbd} \eeq To
estimate the probability that the maximum of $n^2$ excursions is
$> \log n$ we recall (\ref{maxbd}) and that Stirling's formula
implies $k! \ge \delta_0 k^k e^{-k}/\sqrt{2\pi k}$ for some
$\delta_0>0$, so
$$
(\log n)! \ge \delta_1 (\log n)^{\log n} n^{-1} (\log n)^{-1/2} =
\delta_1 n^{\log \log n - 1} (\log n)^{-1/2}
$$
The right-hand size goes to $\infty$ faster than $n^2 \log_2 n$ so
using (\ref{excbd}) we have
$$
P\left( \max_{1\le j \le a\log_2 n} \, \max_{0 \le t \le c}
\zeta^j_t > \log n \right) \to 0
$$
When the last event does not occur we have
$$
\sum_{k=1}^{n^a} \xi^k_t \le a (\log_2 n) \log n = \frac{a}{\log
2} (\log n)^2
$$
Since $a < 2/3 < \log 2 \approx 0.69$, this gives the desired
result.
\end{proof}

\subsection{Completion of the proof of Theorem \ref{supercrit-CLT}}

The cluster queuing system is the first approximation to the
analysis of the dynamics of the supercritical regime. However, it
ignores customer fragmentation and a number of ``bad events" that
we need to consider in order to give a rigorous proof of Theorem
\ref{supercrit-CLT}. Though {\it a priori} one might expect it to
be difficult to take account of corrections of second order, third
order, $\ldots$, and have nightmares about adding up
infinitely many terms, we were pleasantly surprised to see that
the proof could be completed with a few simple estimates.

\mn The first technical problem to confront is to show that the
total amount of mass upstairs stays positive at any given time so
we can apply our time change. This is done in section
\ref{upstairs-mass-section}.

\mn
The more difficult problem is to control the difference between
the CQS and the real system of clusters. To do this, we need a notational scheme to
verify that we have indeed taken care of all of the relevant events. We call clusters
of size larger than $n^a$ {\it large}, those in the CQS (i.e., those that were
generated by a fragmentation of some large cycle), {\it medium},
and non-giant clusters in the random graph {\it small}. Writing
{\it frag} and {\it coag} as shorthand for fragmentation and
coagulation, we have three {\it frag} and six {\it coag} events to handle:

\mn {\it coag(small,small)} is a natural part of the random graph
so these events are not errors. The fragmentation of small
clusters involves $o(n^{1/2})$ clusters and hence does not
significantly alter this process (see {\it frag}(small) and Lemma
\ref{small-frag}).

\mn {\it coag(small,large)} eliminates a small component, but in
the random graph these correspond to the small cluster being
absorbed into the giant component, so this is not an error.

\mn {\it frag(small)} is easy to take care of due to the duality
principle which asserts that finite clusters in the random graph
at time $c>1$ have the same distribution as clusters at time
$c\rho<1$ where $\rho$ is the probability of no percolation. This
allows use to use our subcritical estimates for fragmentation of
small supercritical clusters. More details are given in Lemma
\ref{small-frag}.

\mn {\it coag(large,large)} We do not care about these events
since we do not need to keep track of the number of cycles
upstairs.

\mn
{\it frag(large)} These are the arrivals in the cluster queuing system

\mn {\it coag(medium,large)} are (almost) the
departures in the cluster queuing system. The problem is that the
next three events can cause clusters to gain weight or split into
two.

\mn {\it coag(medium,medium)} are helpful events since they reduce
the number of customers in the CQS. This does make the
fragmentation rate for the new cluster larger than the sum of the
two previous clusters but Lemma \ref{cqs-fragm} will take care of
this. More importantly, it makes the departure rate of the new
cluster larger. This, applied to {\it coag(medium,medium)} and
{\it coag(medium,small)}, shows that the number of medium clusters
is stochastically bounded by the CQS of section \ref{cqs}, and is
the content of Lemma \ref{medium-cqs}.

\mn
{\it coag(medium,small)} eliminates a small component, but in
the random graph these correspond to the small cluster being
absorbed into the giant component.
Again, this also makes the fragmentation rate larger for the
cluster that gained weight but Lemma \ref{cqs-fragm} will take
care of this.

\mn
{\it frag(medium)} is taken care of by Lemma \ref{cqs-fragm}.

\medskip
To complete the proof it remains to prove the three promised lemmas.

\begin{lemma}
\label{medium-cqs} The number of medium clusters is dominated by
that of the CQS. Therefore there are never more than $(\log n)^2$
medium clusters, and never more than $n^a(\log n)^2$ vertices in
medium clusters.
\end{lemma}

\begin{proof} As was just mentioned, the only differences between
the CQS and the medium clusters are generated by events of type
{\it coag(medium,medium)} and {\it coag(small, medium)}. However
both those events do not increase the number of medium clusters,
and both those events make the death rate of the clusters
concerned higher. Hence we can construct the CQS and the medium
clusters process on the same probability space, in such a way that
the {\it total} number of medium clusters is smaller than that of
the CQS.
\end{proof}

\begin{lemma}
\label{cqs-fragm} The expected number of fragmentations of medium
clusters is at most  $O(n^{2a-1}(\log n)^2)$.
\end{lemma}

\begin{proof}
There are never more than $(\log n)^2$ medium clusters. Since
there are at most $n^a$ vertices per medium clusters the total
number of vertices is at most $n^a(\log n)^2$. The rate at which
those fragmenatations happen is thus bounded by
$$
\left(\frac{n^a (\log n)^2}n \right)\frac{n^a}n
$$
so that the expected number of such fragmentations is indeed
$O(n^{2a-1}(\log n)^2)$.
\end{proof}

\begin{lemma}
\label{small-frag} The number of fragmentations of small
components is $o(n^{1/2})$.
\end{lemma}

\begin{proof}
By a now familiar estimate, the expected
number of fragmentations that produce clusters of size smaller
than $n^p$ at times between $n$ and $n+t$ is at most $2n^{p-1}t$.
So we can ignore fragmentations that (a) produce clusters of size
smaller than $n^{0.45}$ before time $cn/2$ and (b) produce clusters of size
smaller than $n^{0.55}$ at times between $n$ and $n + n^{0.9}$.

If $c>1$ the distribution of nongiant components in the random
graph is given by progeny of a Poisson Galton Watson process with
mean $c$ on the event of its extinction. If we let $\rho$ denote
its extinction probability, then the offspring distribution
conditional on extinction is given by
$$
\frac{1}{\rho} e^{-c} \frac{(c\rho)^k}{k!} = e^{-c\rho}
\frac{(c\rho)^k}{k!}
$$
since $\rho = e^{-c(1-\rho)}$. In short, $PGW(c)$ conditioned on
extinction is $PGW(c\rho)$. The last observation implies that
results for finite supercritical clusters can be derived from
those for subcritical clusters. In particular, by Lemma
\ref{tailbound}, the largest nongiant components seen after time $
n + n^{0.9}$, are smaller than $n^{0.2}$. Since fragmentations of
such clusters necessarily produce pieces smaller than $n^{0.2}$
these fragmentations can be ignored by (a).
\end{proof}

\subsection{The initial mass upstairs}
\label{upstairs-mass-section} The last step in the proof of
Theorem \ref{supercrit-CLT} is to ensure that upstairs never
becomes empty in this process. In other words we must prove that
$N_t^{\uparrow}>0$ for all $t>0$ with high probability, so that we
can indeed time-change the queues by $(N_t^{\uparrow})^{-1}$, and
use rigorously all the analysis carried out on (CQS) in section
\ref{cqs}. This will be done by showing that initially there are
already more vertices upstairs than will ever (with high
probability) be taken away by fragmentation in the cluster queuing
system.

\begin{lemma} \label{upstairs-mass-initial}
Initially, upstairs contains at least $N^{\uparrow}_0 \ge
Kn^{1-a/2} $ vertices. In particular $N^{\uparrow}_0 > n^a (\log
n)^2$ and it never becomes empty during the supercritical regime.
\end{lemma}

\mn \begin{proof} Lemma \ref{Stirling} implies that when $c=1$ the
expected number of trees of size $k$
$$
ET_k \sim \frac{n k^{-5/2}}{\sqrt{2\pi}} \exp(-k^3/3n^2)
$$
If we let $|\Cc_{\ge a}| = \sum_{k=n^a}^\infty T_k$ then it
follows that
$$
E|\Cc_{\ge a}| \sim \frac{n}{\sqrt{2\pi}} \sum_{k=n^a}^{\infty}
k^{-5/2} \sim \frac{2}{3\sqrt{2\pi}} n^{1-3a/2}
$$
Bollob\'as (1985) has calculated (see page 107) that the expected
number of ordered pairs of trees of sizes $j$ and $k$,
$$
E(T_j,T_k) \le ET_j ET_k
$$
When $j\neq k$ this implies $\cov(T_j,T_k) \le 0$ and for $j=k$
that $ET_k(T_k-1) \le (ET_k)^2$ or $\var(T_k) \le ET_k$. Summing
we have
$$
\var(|\Cc_{\ge a}|) \le E|\Cc_{\ge a}|
$$
and it follows from Chebyshev's inequality that $|\Cc_{\ge
a}|/E|\Cc_{\ge a}| \to 1$ in probability. These trees have not
experienced fragmentation so their size is always at least $n^a$ and
the total mass in large components is at least $K n^{1-a/2}$. When
$a <2/3$ and $n$ is large, this is much larger than the $n^a(\log
n)^2$ upper bound on the missing mass due to fragmentations.

\mn At this point the proof of Theorem \ref{supercrit-CLT} is
complete.
\end{proof}

\subsection{A sharper estimate for the mass upstairs}

In section \ref{upstairs-mass-section} above, we have just proved
that upstairs never becomes empty in the supercritical regime
(Lemma \ref{upstairs-mass-initial}). But, as was already mentioned
earlier, we expect by analogy with the random graph that in fact a
positive fraction of all $n$ vertices stay upstairs. This is the
content of Theorem \ref{upstairs-mass}, which we restate here for
convenience and then prove.

\mn \textbf{Theorem \ref{upstairs-mass}} \emph{ For any $c>1$, at
time $cn/2$ there are at least $\theta(c) n -o(n)$ vertices
located on large cycles (i.e., of size  greater than or equal to
$n^a$, for any $a<2/3$).}

\begin{proof} In fact it is a simple consequence of Lemmas \ref{medium-cqs}
and \ref{cqs-fragm}. Indeed, the mass missing upstairs must be a
piece of the random graph's giant component fallen downstairs by
fragmentation. Therefore either it is a medium cluster or it has
experienced a consecutive fragmentation. But we now know that
there are never more than $n^a(\log n)^2$ vertices in medium
clusters by Lemma \ref{medium-cqs}. On the other hand, by Lemma
\ref{cqs-fragm}, the expected number of vertices in clusters
having experienced multiple fragmentation has to be smaller than
$$
n^a\cdot Kn^{2a-1}(\log n)^2 = o(n)
$$
as long as $a<2/3$. \end{proof}

\section*{REFERENCES}

\mn Aldous, D. (1997) Brownian excursions, critical random graphs
and the multiplicative coalescent. {\it Ann. Prob.} 25, 812--854

\mn Aldous, D. (1999) Deterministic and stochastic models for
coalescence (aggregation and coagulation) : a review of the
mean-field theory for probabilists. {\it Bernoulli.} 5, 3--48

\mn Angel, O. (2003) Random infinite permutations and the cyclic
time random walk. Pages 9--16 in Banderier and Krattenthaler
(2003)

\mn Arratia, R. and Barbour, A. and Tavar\'e, S. (2003) {\it
Logarithmic combinatorial structures : a probabilistic approach.}
European Math. Society Monographs, 1.

\mn Bafna, V. and Pevzner, P. (1995) Sorting by reversals: Genome
rearrangement in plant organelles and evolutionary history of X
chromosome. {\it Mol. Biol. Evol.} 12, 239--246

\mn Banderier, C., and Krattenthaler, C. (2003) Proceedings of the
conference Discrete Random Walks. {\it Discrete Math and Computer
Science.} \texttt{dmtcs.loria.fr/proceedings/dmACind.html}

\mn Berestycki, N. and Durrett, R., (2003) A phase transition in
the random transposition random walk. Pages 17-26 in Banderier and
Krattenthaler (2003)

\mn Bollob\'as, B. (1984) The evolution of random graphs. {\it
Trans. Amer. Math. Soc.} 286, 257--274

\mn Bollob\'as, B. (1985) {\it Random Graphs}, Cambridge
University Press.

\mn Borel, E. (1942) Sur l'emploi du th\'eor\`eme de {B}ernoulli
pour faciliter le calcul d'une infinit\'e de coefficients.
Application au probl\`eme de l'attente \`a un guichet. {\it C.R.
Acad. Sci. Paris.} 214, 452--456

\mn Bourque, G. and Pevzner, P. A. (2002) Genome-scale evolution:
reconstructing gene orders in the ancestral species. {\it Genome
Research}. 12, 26--36

\mn Devroye, L. (1992) The branching process method in the
Lagrange random variate generation,
\texttt{cgm.cs.mcgill.ca/\~{}luc/branchingpaper.ps}

\mn Diaconis, P., Mayer-Wolf, E., Zeitouni, O., and Zerner, M.
(2003) Uniqueness of invariant distributions for split-merge
transformations and the Poisson-Dirichlet law. {\it Ann. Prob.},
to appear

\mn Durrett, R. (1996) {\it Probability: Theory and Examples},
Second Edition, Duxbury Press

\mn Durrett, R. (2002) {\it Probability Models for DNA Sequence
Evolution.} Springer-Verlag, New York

\mn Durrett, R. (2003) Shuffling Chromosomes. {\it J. Theor.
Prob.} 16, 725--750

\mn Durrett, R., Nielsen, R., and York, T.L. (2003) Bayesian
estimation of genomic distance. {\it Genetics}, to appear

\mn Hannehalli, S. and Pevzner, P.A. (1995) Transforming cabbage
into turnip (polynomial algorithm for sorting signed permutations
by reversals). {\it Proceedings of the $27^{th}$ Annual Symposium
on the Theory of Computing}, 178--189. Full version in the
\emph{Journal of the ACM.} 46, 1--27

\mn Jacod, J. and Shiryaev, A. (1987) {\it Limit Theorems for
Stochastic Processes}, Springer New-York

\mn Janson, S., Knuth, D. E., Luczak, T. and Pittel, B. (1993) The
birth of the giant component. {\it Rand. Struct. Algor.} 4,
231--358

\mn Janson, S., Luczak, T., and Ruczinski, A. (2000) {\it Random
Graphs}, Wiley-Interscience, New York

\mn Luczak, T., Pittel, B., and Wierman, J. C. (1994) The
structure of a random graph near the point of the phase
transition. {\it Trans. Amer. Math. Soc.} 341, 721--748

\mn Mayer-Wolf, E. and Zeitouni, O. and Zerner, M. (2002)
Asymptotics of certain coagulation-fragmentation processes and
invariant Poisson-Dirichlet measures. {\it Electr. Journ. Prob.}
7, 1--25

\mn Pevzner, P.A. (2000) {\it Computational Molecular Biology: An
Algorithmic Approach.} MIT Press, Cambridge

\mn Pevzner, P.A. and Tesler, G. (2003) Genome rearrangement in
mammalian evolution: lessons from human and mouse genomes. {\it
Genome Research}. 13, 37--45

\mn Pitman, J. (1998) Enumerations of trees and forests related to
branching processes and random walks. \emph{Microsurveys in
Discrete Probability}, D. Aldous and J. Propp editors. DIMACS Ser.
Discrete Math. Theoret. Comp. Sci no.41 163-180. Amer. Math. Soc.
Providence RI.

\mn Pitman, J. (1999) Coalescent random forests, \emph{J. Comb.
Theory A}. 85 165-193.

\mn Pitman, J. (2002) Poisson-Dirichlet and GEM invariant
distributions for split-and-merge transformations of an interval
partition.  {\it Combin. Prob. Comput.} 11, 501--514

\mn Pitman, J. (2003) Combinatorial stochastic processes. {\it
Lecture Notes for St. Flour Course}. To appear, available at
{\texttt{http://stat-www.berkeley.edu/users/pitman/}}

\mn Pittel, B. (1990) On tree census and the giant component in
sparse random graphs, {\it Rand. Struct. Algor.}, {\bf 1},
311--342

\mn Ranz, J.M. and Casals, F. and Ruiz, A. (2001) How malleable is
the eukaryotic genome? Extreme rate of chromosomal rearrangement
in the genus \emph{Drosophila}. {\it Genome Research}. 11,
230--239

\mn Revuz, D. and Yor, M., (1999) {\it Continuous martingales and
Brownian Motion}, Springer-Verlag, New York

\mn Schramm, O. (2004) Composition of random transpositions, to
appear.

\mn Tanner, J.C. (1961) A derivation of the Borel distribution.
{\it Biometrika} 48, 222--224

\mn York, T.L., Durrett, R., and Nielsen, R. (2002) Bayesian
estimation of inversions in the history of two chromosomes. {\it
J. Comp. Bio.} 9,808--818

\clearpage

\begin{figure}
\label{breakpoint}
\begin{center}
\begin{picture}(370,70)
\put(13,5){\erm 0} \put(26,5){\erm 1} \put(43,5){\erm 2}
\put(53,5){\erm 14} \put(70,5){\erm 13} \put(85,5){\erm 11}
\put(99,5){\erm 12} \put(115,5){\erm 20} \put(130,5){\erm 19}
\put(145,5){\erm 17} \put(160,5){\erm 18} \put(175,5){\erm 16}
\put(190,5){\erm 15} \put(208,5){\erm 3} \put(223,5){\erm 4}
\put(235,5){\erm 22} \put(251,5){\erm 21} \put(266,5){\erm 6}
\put(283,5){\erm 5} \put(296,5){\erm 9} \put(311,5){\erm 10}
\put(327,5){\erm 7} \put(343,5){\erm 8} \put(354,5){\erm 23}
\put(15,15){\line(0,1){10}} \put(30,15){\line(0,1){10}}
\put(15,25){\line(1,0){15}} \put(225,15){\line(0,1){10}}
\put(285,15){\line(0,1){10}} \put(225,25){\line(1,0){60}}
\put(300,15){\line(0,1){10}} \put(345,15){\line(0,1){10}}
\put(300,25){\line(1,0){45}} \put(75,15){\line(0,1){10}}
\put(105,15){\line(0,1){10}} \put(75,25){\line(1,0){30}}
\put(150,15){\line(0,1){10}} \put(180,15){\line(0,1){10}}
\put(150,25){\line(1,0){30}} \put(135,15){\line(0,1){15}}
\put(165,15){\line(0,1){15}} \put(135,30){\line(1,0){30}}
\put(270,15){\line(0,1){15}} \put(330,15){\line(0,1){15}}
\put(270,30){\line(1,0){60}} \put(60,15){\line(0,1){20}}
\put(195,15){\line(0,1){20}} \put(60,35){\line(1,0){135}}
\put(45,15){\line(0,1){25}} \put(210,15){\line(0,1){25}}
\put(45,40){\line(1,0){165}} \put(240,15){\line(0,1){20}}
\put(360,15){\line(0,1){20}} \put(240,35){\line(1,0){120}}
\put(120,15){\line(0,1){30}} \put(255,15){\line(0,1){30}}
\put(120,45){\line(1,0){135}} \put(90,15){\line(0,1){35}}
\put(315,15){\line(0,1){35}} \put(90,50){\line(1,0){225}}
\linethickness{0.7mm} \put(15,15){\line(1,0){15}}
\put(45,15){\line(1,0){15}} \put(75,15){\line(1,0){15}}
\put(105,15){\line(1,0){15}} \put(135,15){\line(1,0){15}}
\put(165,15){\line(1,0){15}} \put(195,15){\line(1,0){15}}
\put(225,15){\line(1,0){15}} \put(255,15){\line(1,0){15}}
\put(285,15){\line(1,0){15}} \put(315,15){\line(1,0){15}}
\put(345,15){\line(1,0){15}}
\end{picture}
\caption{Breakpoint graph for human-mouse X chromosome comparison}
\end{center}
\end{figure}
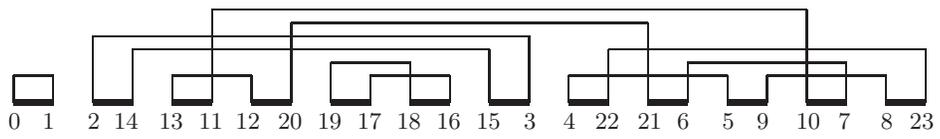

\clearpage

\begin{table}[tp] \label{second_data_set}
\begin{center}
\begin{tabular}{cccccccccc}
36 & 37 & 17 & 40 & 16 & 15 & 14 & 63 & 10 & 9 \cr 55 & 28 & 13 &
51 & 22 & 79 & 39 & 70 & 66 & 5 \cr 6 &  7 & 35 & 64 & 33 & 32 &
60 & 61 & 18 & 65 \cr 62 & 12 &  1 & 11 & 23 & 20 & 4 & 52 & 68 &
29 \cr 48 &  3 & 21 & 53 &  8 & 43 & 72 & 58 & 57 & 56 \cr 19 & 49
& 34 & 59 & 30 & 77 & 31 & 67 & 44 & 2 \cr 27 & 38 & 50 & 26 & 25
& 76 & 69 & 41 & 24 & 75 \cr 71 & 78 & 73 & 47 & 54 & 45 & 74 & 42
& 46
\end{tabular}
\end{center}
\caption{Order of the genes in {\it D. repleta} compared to their
order in {\it D. melanogaster}}
\end{table}

\end{document}